\documentclass{amsart}
\usepackage{amsmath,amsthm,amsfonts,amssymb}
\usepackage{float}

\usepackage{graphicx}

\floatstyle{ruled}\restylefloat{figure}

\usepackage{psfrag}
\newcommand{\mkfrag}[1]{\psfrag{#1}{$#1$}}
\mkfrag{t}
\mkfrag{u}\mkfrag{u'}
\mkfrag{v}
\mkfrag{w}
\mkfrag{z}
\mkfrag{a}\mkfrag{a'}
\mkfrag{b}\mkfrag{b'}
\mkfrag{c}

\newcommand{\Ints}{\mathbb{Z}}
\newcommand{\Nats}{\mathbb{N}}
\newcommand{\R}{\mathbb{R}}
\newcommand{\Expectation}{\mathbb{E}}
\newcommand{\Expectationp}[2]{\mathbb{E}_{#1} {#2}}
\newcommand{\Prob}[1]{\mathbb{P}\{#1\}}
\newcommand{\Probp}[2]{\mathbb{P}_{#1}\{#2\}}

\newcommand{\aut}{\operatorname{Aut}}
\newcommand{\supp}{\operatorname{supp}}
\newcommand{\tree}[1]{#1}
\newcommand{\st}{:}
\newcommand{\defn}[1]{\emph{#1}}
\newcommand{\bydef}{:=}

\newcommand{\bdy}[1]{\partial_E #1}
\newcommand{\stab}[1]{\textrm{Stab}({#1})}
\newcommand{\symdiff}{\mathbin\vartriangle}

\newcommand{\A}{\mathcal{A}}
\newcommand{\B}{\mathcal{B}}
\newcommand{\GG}{\mathcal{G}}
\newcommand{\FF}{\mathcal{F}}

\newcommand{\dist}{\operatorname{dist}}

\newcommand{\chap}{\section}
\newcommand{\sect}{\subsection}

\newcommand{\sizeof}[1]{{\left|#1\right|}}

\newcommand{\set}[1]{{\left\{#1\right\}}}

\newtheorem{theorem}{Theorem}[section]
\newtheorem{conjecture}[theorem]{Conjecture}
\newtheorem{lemma}[theorem]{Lemma}
\newtheorem{xampl}[theorem]{Example}
\newtheorem{note}[theorem]{Note}
\newtheorem{definition}[theorem]{Definition}

\begin{document}

\title{Critical Percolation on Certain Non-unimodular Graphs}

\author{Yuval Peres} \email{peres@stat.berkeley.edu\qquad
http://www.stat.berkeley.edu/\~{}peres} \address{Departments of Statistics
and Mathematics, 367 Evans Hall, University of California, Berkeley, CA
94720}

\author{G\'abor Pete}
\email{gabor@stat.berkeley.edu\qquad http://www.stat.berkeley.edu/\~{}gabor}
\address{Department of Statistics, 367 Evans Hall, University of
California, Berkeley, CA 94720}

\author{Ariel Scolnicov}
\email{ascolnic@checkpoint.com}

\date{October 25, 2005}

\thanks{Our work was partially supported by the NSF grant DMS-0244479
(Peres, Pete), and OTKA (Hungarian National Foundation for Scientific
Research) grants T30074 and T049398 (Pete).}

\begin{abstract}
 An important conjecture in percolation theory is that almost surely no
 infinite cluster exists in critical percolation on any transitive graph for
 which the critical probability is less than 1.  Earlier work has established
 this for the amenable cases $\Ints^2$ and $\Ints^d$ for large~$d$, as well
 as for all non-amenable graphs with unimodular automorphism groups.
 We show that the conjecture holds for the basic classes of non-amenable
 graphs with non-unimodular automorphism groups: for decorated trees and
 the non-unimodular
 Diestel-Leader graphs. We also show that the connection probability 
 between two vertices
 decays exponentially in their distance. Finally, we prove that critical 
 percolation on the
 positive part of the lamplighter group has no infinite clusters. 
\end{abstract}

\maketitle

\chap{Introduction and Preliminaries}
\label{chap:intro}

\sect{Introduction}
\label{sec:intro:perc}

We will focus on the following general conjecture of
Benjamini and Schramm \cite{BS:beyond-Z^d:1996} on critical
percolation (see Subsection \ref{sec:back} for definitions):

\begin{conjecture}
\label{conj:intro:transitive}
Let~$G$ be a transitive graph.  If $p_c<1$, then almost surely
critical percolation on~$G$ has no infinite clusters.
\end{conjecture}

Earlier work of Harris \cite{Harris:1960} and Kesten \cite{Kesten:1980}
established that for the graph $\Ints^2$ critical percolation almost
surely has no infinite cluster at $p=p_c$.  Later, Hara and Slade
\cite{HaSl:1994} established the same for $\Ints^{d}$, when $d\ge 19$.
However, Conjecture \ref{conj:intro:transitive} remains open
for~$\Ints^{d}$ where $3\le d\le 18$, along with many other amenable
graphs.

For regular trees, the conjecture is just the classical result that a
critical Galton-Watson tree dies out. Wu \cite{Wu:1993} showed the 
conjecture for products of regular trees and $\Ints$, when the degree 
of the regular tree was large enough.

Benjamini, Lyons, Peres and Schramm, see \cite{BLPS:1999} or
\cite{BLPS:gafa:1999}, proved the conjecture for non-amenable graphs when
the automorphism group of the graph is unimodular. However, their ``Mass
Transport Principle'' does not adapt to the case when the action is not
unimodular, leaving that case open, as well.

In this paper, we prove the conjecture for the best-known examples of
transitive graphs with non-unimodular automorphism groups. 

\begin{figure}[htbp]
  \begin{center}
    \includegraphics{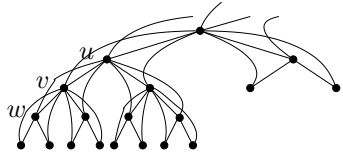}
    \caption{``Grandparent'' tree.}
    \label{fig:grandfather}
  \end{center}
\end{figure}

Section \ref{chap:dt} deals with trees ``decorated'' by adding
edges, where this decoration retains transitivity.  A typical example is 
the ``grandparent tree'' (see Figure \ref{fig:grandfather}),  
due to Trofimov \cite{Tr:1985}, and also appearing
in Soardi and Woess \cite{SoWo:1990}.

Section \ref{chap:43} proves the conjecture for a class of graphs due to
Diestel and Leader \cite{DiLe:non-cayley}, see Figure \ref{fig:43}. Such a
graph $\Gamma_{\alpha,\beta}$ is the ``horocyclic product'' of an
$(\alpha+1)$-regular tree $T_\alpha$ and a $(\beta+1)$-regular tree
$T_\beta$; see Subsection \ref{sec:43:defn} for a formal definition.

When $\alpha\ne\beta$, the graph~$\Gamma_{\alpha,\beta}$ is
non-unimodular. Theorem \ref{thm:43:dead} proves Conjecture
\ref{conj:intro:transitive} for these graphs $\Gamma_{\alpha,\beta}$.

If $\alpha=\beta$, the graph $\Gamma_{\alpha,\beta}$ turns out to be a
Cayley graph of the ``lamplighter group'' of Ka{\u\i}manovich and Vershik
(Example 6.1 of \cite{KV:lamplighter:1983}).  As such, it is unimodular,
moreover, it is amenable, which means that we are unable to prove the
conjecture in this case.  However, in Section \ref{chap:lamp} we show that
critical percolation on the positive part of this graph, and also of
another natural Cayley graph of the lamplighter group, has no infinite
components. This is analogous to the case of half-space percolation in
$\Ints^d$, see \cite{BGN:half-space:1991}.

We also show that, in critical percolation on any of our non-unimodular
graphs, the connection probability between two vertices decays
exponentially in their distance. The importance of such an exponential
decay is discussed in \cite{BS:recent}; for example, it might help in
proving the existence of the non-uniqueness phase. For the Diestel-Leader
graphs $\Gamma_{\alpha,\beta}$, the standard methods give the existence of
this phase only when $\beta$ is sufficiently large compared to $\alpha$
(or vice versa), see at the end of Section~\ref{chap:43}.

Note that the method of \cite{ALy:1991}, see also \cite[Corollary
5.5]{BLPS:gafa:1999}, shows that if the automorphism group of a transitive
graph is non-amenable (which is stronger than the non-amenability of the
graph), then there can not be a unique infinite cluster in critical
percolation. However, our examples have amenable automorphism groups,
hence our proofs have to rule out the possibility of a unique infinite
cluster, as well, which is non-trivial in the case of the Diestel-Leader
graphs. A very recent preprint of \'Ad\'am Tim\'ar \cite{Adam:nonuni},
together with an unpublished result of Lyons, Peres and Schramm, 
show that there cannot exist infinitely many infinite clusters in
critical percolation on any non-unimodular transitive graph. 

\sect{Background: amenability, unimodularity, and percolation}
\label{sec:back}

Let $G=(V,E)$ be a locally finite infinite graph. Denote by
$\aut{G}$ its group of automorphisms, i.e.~the group of bijective 
maps $g:V(G)\longrightarrow V(G)$ such that $\{u,v\}\in E(G)$ 
iff $\{gu,gv\}\in E(G)$. $G$ is called \defn{transitive} if any 
pair of vertices of~$G$ has an automorphism that maps the first vertex 
to the second one.

If we equip $\Gamma=\aut{G}$ with the topology of pointwise convergence on
$G$, then it becomes a locally compact topological group.  Therefore
it has both left- and right-invariant Haar measures, and we can
consider the Banach space $L^{\infty}(\Gamma)$ of measurable essentially
bounded real valued functions on $\Gamma$ w.r.t.~the left-invariant Haar
measure. A linear functional $m:\,L^{\infty}(\Gamma)\longrightarrow \R$ is
called an {\it invariant mean} if it maps nonnegative functions to
nonnegative reals, the constant ${\mathbf 1}$ function to 1, and
$m(L_g\phi)=m(\phi)$ for any $g\in \Gamma$ and $\phi\in L^{\infty}(\Gamma)$,
where $L_g(\phi)(h):=\phi(gh)$.

\begin{definition}\hfil
  \begin{itemize}
  \item The \defn{edge-isoperimetric constant} of a graph $G$ is
    \[
    \iota_E(G) \bydef \inf \set{\frac{\sizeof{\bdy K}}{\sizeof{K}}\st
      K\subset V(G), \sizeof{K} < \infty},
    \]
    where $\bdy K \bydef \set{\set{u,v}\in E(G)\st u\in K, v\notin K}$.
$G$ is \defn{amenable} if $\iota_E(G) = 0$.
  \item A locally compact topological group $\Gamma$ is \defn{amenable} if
it has an invariant mean. If $\Gamma$ is finitely generated, then this is
equivalent to saying that it has an amenable Cayley graph, see
\cite{amenbook}.
   \item A locally compact topological group is called 
\defn{unimodular} if its left- and right-invariant Haar-measures coincide.
  \end{itemize}
\end{definition}

Schlichting \cite{Sc:1979} and Trofimov \cite{Tr:1985} give a
combinatorial characterization of unimodularity, which is made
explicit by Soardi and Woess \cite{SoWo:1990} for the action of a
group of graph automorphisms on the graph. According to this
characterization, the action of a group of automorphisms~$\Gamma$ on a
graph~$G$ is unimodular if and only if for any pair $x,y$ of
vertices, $\sizeof{\stab{x}\cdot y} = \sizeof{\stab{y}\cdot x}$, where
$\stab{x} = \set{g\in\Gamma \st gx = x}$ is the \defn{stabilizer}
of~$x$. We say that a transitive graph $G$ is \defn{unimodular} if the
action of the full group $\aut G$ is unimodular.

There are basic connections between non-amenability of a graph and the
non-amenabi\-li\-ty of its automorphism group.  A useful lemma, see
e.g.~Lemma 3.3 of \cite{BLPS:gafa:1999}, allows us to take invariant 
means on \emph{any} appropriate graph, instead of on the group
itself. If $G$ is a countable graph, and $\Gamma$ is a closed subgroup
of $\aut G$, then $\Gamma$ acts on the Banach space $L^\infty(V(G))$ of
real valued bounded functions by $L_g(\phi)(v):=\phi(gv)$, and we
can define a $\Gamma$-invariant mean on $G$ analogously to how we did
above. 

\begin{lemma}[Characterization of group-amenability]
\label{lem:amenchar}
  Let $G$ be a graph and $\Gamma$ be a closed subgroup of $\aut G$.
 Then $\Gamma$ is amenable if and only if $G$ has a $\Gamma$-invariant
mean.\qed 
\end{lemma}

There is also a characterization of graph amenability in terms of the
amenability of closed transitive groups of automorphisms, due to
Soardi and Woess \cite{SoWo:1990}. See also Theorem 3.4 of 
\cite{BLPS:gafa:1999}.

\begin{lemma}[Corollary 1 of \cite{SoWo:1990}]
\label{lem:amenunim}
  Let~$G$ be a graph and~$\Gamma$ a closed transitive subgroup
  of $\aut{G}$.  Then $G$ is amenable if and only if $\Gamma$ is amenable
and its action is unimodular.\qed 
\end{lemma}

Given a graph $G$ and $0\le p\le 1$, \defn{percolation} on~$G$ is a
measure $\Probp{p}{\cdot}$ on subsets $\mathcal{E}\subseteq E(G)$, where
the events $\set{e\in\mathcal{E}}$, $e\in E(G)$, are all independent and
occur with probability~$p$.  Edges $e\in\mathcal{E}$ are called
\defn{open}, and edges $e\in\mathcal{E}^{\textrm{c}}$ \defn{closed}; paths
shall be called \defn{open} if all edges are open.  The \defn{cluster} of
a vertex $o\in V(G)$ is 
\[ C(o) = \set{v\st \text{$o\longleftrightarrow v$
by an open path}}. \] 
By Kolmogorov's 0-1 law, for any value of~$p$ an
infinite cluster exists with probability 0 or 1. So, define the
\defn{critical probability} $p_c$ for percolation by \[ p_c = \inf
\set{p\st \Probp{p}{\text{$\exists\ \infty$ cluster}} = 1}. \]

When the value of $p$ is clear from the context, and especially when
$p=p_c$, we write $\Prob{\cdot}$ for $\Probp{p}{\cdot}$.

For non-amenable graphs with bounded degree it is known
\cite{BS:beyond-Z^d:1996} that $0<p_c<1$, hence Conjecture
\ref{conj:intro:transitive} poses a real question in this case.

For non-amenable transitive graphs there is a second critical value of
interest, \[ p_u= \inf \set{p\st \Probp{p}{\exists\text{ a unique $\infty$
cluster}} = 1}. \] Another famous conjecture of \cite{BS:beyond-Z^d:1996}
is the strict inequality $p_c<p_u$ for these graphs. The standard
references for percolation are \cite{grimmbook} and \cite{LPbook}.

\sect{The general strategy}
\label{sec:strategy}

The main steps of the proof are shared by all the examples we deal
with.  First, we shall use the tree structure underlying the graph to
construct a Galton-Watson process and bound the 
expected number of vertices at level $k$ that can be reached
from a fixed vertex $o$ at level 0 via certain restricted paths
that stay in the ``downwards half-graph'' from $o$.
 Then a Fatou lemma
argument will imply that the component is a.s.~finite in this
downwards half-graph. Moreover, as the combinatorial characterization
of non-unimodularity suggests, the component of a vertex has more ways
to grow ``downwards'' than ``upwards'', so the component cannot
directly reach infinitely far upwards, either. In a decorated tree
there is no ``sideways'' direction, so it follows easily that the
entire component must be finite. For the Diestel-Leader graphs the
specific combinatorial structure helps in showing that the
``exponentially unlikely'' upward growth makes it impossible that
there is a cluster oscillating infinitely up and down.

\chap{Decorated trees}
\label{chap:dt}

\sect{Definition and examples}
\label{sec:dt:defn}

Let~$\tree{T}$ be a $d+1$-regular tree.  $\tree{T}$~is a transitive
non-amenable graph, $\aut \tree{T}$~is non-amenable, and its action
on~$\tree{T}$ is unimodular.  We shall examine a class of non-amenable
transitive graphs~$G$ derived from~$\tree{T}$ by adding edges to it
for which $\aut G$ will be amenable (and therefore, by Lemma
\ref{lem:amenunim}, will act on $G$ in a non-unimodular manner).  

Two rays (half-infinite simple paths) in $T$ are called equivalent if they
differ only in finitely many edges. An {\it end} of the tree is an
equivalence class of rays. Pick an end $\xi$ of $\tree{T}$ and direct all
edges of $\tree{T}$ towards~$\xi$.  If there is an edge from~$v$ to~$u$,
we say that~$v$ is the \defn{child} of~$u$, and~$u$ is its \defn{parent}.  
We shall use the terms \defn{sibling}, \defn{grandchild} and
\defn{grandparent} in their obvious meaning.  We say that~$v$ is a
\defn{descendant} of~$u$ and that~$u$ is an \defn{ancestor} of~$v$ if
there is a directed path from~$v$ to~$u$.  The \defn{downwards
subtree}~$S_v$ of a vertex~$v$ is the graph on the vertices descended
from~$v$.  Distinguishing some vertex~$o\in\tree{T}$, we may define a
\defn{level function} $\ell: V(\tree{T})\to \Ints$ by $\ell(o)=0$ and
$\ell(v)=\ell(u)+1$ whenever $v$ is a child of $u$. Note that large values
of this level function mean large depths in $T$, while negative values
correspond to being higher than $o$. When considering the cluster of a
given vertex~$o$, we shall frequently make use of the \defn{level sets}
(relative to~$o$), defined for $k\in\Ints$ by $L_k \bydef \{ v \st
\ell(v)-\ell(o)=k \}$. For example, the visually clear expression that a
path $v_1,\ldots,v_n$ \defn{does not go above level $L_k$} can be written
as $\ell(v_i)\geq k$ for all $0\leq i\leq n$.

Let $K=\{\alpha\in\aut \tree{T} \st \alpha\xi=\xi\}$ be the group of
$\xi$-preserving automorphisms of~$\tree{T}$.  Then~$K$ is an
\emph{amenable} group (any Banach limit on~$\xi$ is a 
$K$-invariant mean, which suffices by Lemma \ref{lem:amenchar}),  
which acts on~$\tree{T}$ transitively.

Now let~$L$ be some subgroup of~$K$ (possibly~$K$ itself) which acts
transitively on~$\tree{T}$.  Any locally finite graph $G=(V(\tree{T}),
E(\tree{T})\cup E')$ with $L\cdot E' = E'$ will be called a
\defn{decorated tree} (or \defn{$L$-decorated tree}).  The graph~$G$
itself is always non-amenable, since it results from the non-amenable
graph~$\tree{T}$ by adding edges.  Considering the action of~$L$ on
the vertices of $G$, we may regard it as a subgroup of~$\aut G$;
however, $\aut G$ might still be non-amenable.

\begin{figure}[htbp]
  \begin{center}
    \includegraphics{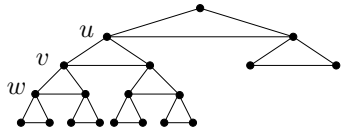}
    \caption{``Triangles'' tree.}
    \label{fig:triangle}
  \end{center}
\end{figure}

\begin{xampl}[$\aut G$ non-amenable, unimodular action]
\label{ex:dt:triangle}
Take $d=2$ and $L=K$, and let 
$E' = \set{\set{u,v} \st \text{$u$,$v$ are siblings}}$, see Figure 
\ref{fig:triangle}.  Then~$\aut G$ is non-amenable, 
and its action is unimodular.
\end{xampl}

\begin{xampl}[$\aut G$ amenable, non-unimodular action]
\label{ex:dt:grandfather}
  Take $L=K$, and $E'=\{\{u,w\} \st u$ is the grandparent of
  $w\}$, see Figure \ref{fig:grandfather}. The action of $\aut
  G$ on~$G$ is not unimodular, and it is an amenable group.
\end{xampl}

For the remainder of this section, we shall fix some graph~$G$ which
is a decorated tree, and prove that critical percolation on~$G$ almost
surely has no infinite components.  While all the results hold for
$\aut G$ with unimodular action, this case is covered by 
\cite{BLPS:gafa:1999}; the result is new only for $\aut G$
with non-unimodular action.

\begin{theorem}\label{thm:dt:dead}
  Let~$G$ be a decorated tree.  Then critical percolation on~$G$ 
  a.s.~has no infinite components. 
\end{theorem}

\begin{note} \rm It is easy to check directly, and also
follows from our proof of Theorem \ref {thm:dt:dead} below,
 that for any $p<1$, percolation on a decorated tree satisfies
 $\Probp{p}{x\longleftrightarrow y}\leq e^{-c \dist (x,y)}$ for some $c=c_p>0$. 
The same
result at $p=p_c$ for the Diestel-Leader graphs is not that easy, 
and will be proven in Section \ref{chap:43}.
\end{note}

We will follow the strategy outlined in Subsection~\ref{sec:strategy}.

\sect{Bounding branches}
\label{sec:dt:branch_bound}

\begin{definition}
\label{def:dt:forward}
  Consider percolation on a decorated tree $G$. The \defn{forward cluster}
$C^{+}(o)$ of a vertex $o\in V(G)$ is defined by
  \[
\begin{split}
  C^{+}(o) \bydef \{v\st v\longleftrightarrow o 
& \text{ by an open path }(o=v_0,v_1,\dots,v_n=v)\\ 
&\text{ inside the downwards subtree }S_o,\\
&\text{ with }\ell(v_i)\leq\ell(v)\text{ for all }0\leq i< n\}.
\end{split}
  \]
\end{definition}

Note that $C^+(o)$ is not necessarily connected. We start by showing that
$C^+(o)$ is ``narrow'', in the sense that it contains few branches.

\begin{lemma}\label{lem:dt:bb}
  Consider critical percolation on a decorated tree $G$ at $p_c(G)$. Let
$o\in V(G)$ be a vertex, and
  define $V^+_k \bydef C^+(o) \cap L_k $.
  Then $e_k \bydef \Expectation{\sizeof{V^+_k}} \le 1$ for all~$k \ge 0$.
\end{lemma}

\begin{proof}
  Suppose to the contrary that~$e_{k_0} > 1$ for some $k_0\ge 0$.  We shall
  use this to find subsets of the vertices of $\{V^+_{j\cdot
  k_0}\}_{j=0}^{\infty}$ which will form a supercritical Galton-Watson
  process:

  \begin{itemize}
  \item The root of the process shall be the vertex~$Y_{0,1} = \{o\}$.
  \item If at level~$j-1$ of the process we picked vertices
    $v_{j-1,1},\ldots,v_{j-1,N_{j-1}} \in V^+_{(j-1)k_0}$, we shall
    pick at level $j$ as descendants of each $v_{j-1,i}$ the vertices
    $Y_{j,i} = C^+(v_{j-1,i}) \cap L_{j k_0}$.
  \end{itemize}
  
  Due to the construction, for any fixed~$j$, if we condition on the
previous generation $\{Y_{j-1,i}:i=1,\dots,N_{j-1}\}$, then the sets
$Y_{j,i}$ are independent. Also, $\sizeof{Y_{j,i}}$ has the same
distribution as $\sizeof{Y_{1,1}}$, so this is indeed a Galton-Watson
process.
  
  Since $Y_{1,1} = V^+_{k_0}$, this is a supercritical process.  But
  $e_{k_0}$ is a polynomial in~$p$, and in particular is continuous.
  Thus, we may decrease~$p$ below $p_c$ keeping $e_{k_0}(p) > 1$.
  This would give a positive probability for
  \[
  \sizeof{C(o)} \ge
  \sum_{j=0}^{\infty} \sizeof{V^+_{j\cdot k_0}} = \infty,
  \]
  contradicting criticality at~$p_c$.
\end{proof}

\sect{Clusters are finite}
\label{sec:dt:finite}

Define $r$ as the maximal length of a path in~$\tree{T}$ connecting
the two endpoints of an edge of~$G$.  Since $G$ is locally finite and 
transitive, $r$ is well-defined and finite. Furthermore, let $D$ be the common 
degree of the vertices of $G$.  

\begin{lemma}\label{lem:dt:posfinite}
  In critical percolation on a decorated tree $G$, for any $o\in V(G)$, 
the forward cluster $C^+(o)$ is a.s.~finite.
\end{lemma}

\begin{proof} Consider the band of levels $H_k:=\cup_{j=k}^{k+r} L_j$.
Recall the random variables $|V_k^+|$ from Lemma \ref{lem:dt:bb}. For a
percolation configuration $\omega$, let $E_k$ be the set of edges in open
paths leading from $o$ to $L_k$, staying in $S_o$ and not going below
$L_k$. Define $F_k:=\cup_{j=k}^{k+r} E_j$ and $W^+_k:=\cup_{j=k}^{k+r}
V^+_j$. Then Lemma \ref{lem:dt:bb} and Fatou's lemma give us that 
\[
\Expectation{\bigl\{\liminf_{k\to\infty} |W_k^+|\bigr\}}\leq
\liminf_{k\to\infty}\Expectation{|W_k^+|}\leq r+1. 
\] 
Thus the random
variable $\liminf_{k\to\infty} |W_k^+|$ is almost surely finite.

The event $\A_k(n):=\{|W_k^+|=n\}$ is clearly determined by $F_k$.
Furthermore, given $F_k$ and $\A_k(n)$ with some $n\geq 1$, the
probability of the event $\GG_k:=\{$all edges incident to $W_k^+$ and not
in $F_k$ are closed$\}$ is at least $(1-p_c)^{Dn}>0$. This means that if
$\A_k(n)$ happens for infinitely many $k$ values, then there is almost
surely a $K$ such that $\GG_K$ occurs. But note that if $u\in V_m^+$ for
some $m>k+r$, then any open simple path from $o$ to $u$ that shows this,
when it first enters $H_k$, goes through some vertex $v\in W_k^+$, and it
leaves $W_k^+$ the last time through an edge not in $F_k$. Hence, $\GG_K$
implies that $V_m^+=\emptyset$ for all $m>K+r$, which means that $\A_k(n)$
could not happen infinitely often.

Therefore, we must have $|W_k^+|=0$ infinitely often a.s. But the above
argument also shows that $W_k^+=\emptyset$ implies $V_m^+=\emptyset$ for
all $m>k+r$, hence $|C^+(o)|<\infty$ a.s. 
\end{proof}

In the case of a decorated tree it is particularly easy to use the tree-like 
structure of $G$ to show that clusters cannot extend infinitely far ``up'' or 
``sideways''.

\begin{lemma}\label{lem:dt:subtree} In independent $p$-percolation with
any $p<1$, the cluster $C(o)$ is a.s.~contained in some downwards subtree.
\end{lemma}

\begin{proof}
  Call the subtree $S_v$ of any vertex~$v$ \defn{isolated} if no open edges
  remain connecting $S_v$ with~$V(G)\setminus S_v$; define the events 
  \(I_v =  \{\text{$S_v$ is isolated}\}\).
  
  Recall the bound $r$ on the ``maximal length in~$\tree{T}$'' of an
  edge of~$G$. It follows that $I_v$ depends only on a constant finite
  number of edges. Consider now the events $I_{v_1}, I_{v_2}, \ldots$
  for some vertices $v_1, v_2,\ldots$ on the path upwards from~$o$,
  which are sufficiently far apart so that these events are all
  independent. Then a.s.~one (in fact, infinitely many) of the
  $I_{v_i}$ will occur, and $C(o)$ is contained in the downwards subtree
  of this $v_i$. The probability that the distance of $v_i$ from $o$ is
  larger than $t$ decays exponentially in $t$.  
\end{proof}

\begin{proof}[Proof of Theorem \ref{thm:dt:dead}]
  Almost surely, the conclusions of Lemmas \ref{lem:dt:posfinite} and
  \ref{lem:dt:subtree} hold for all vertices of~$G$. Similarly, it is
  enough to show that $C(o)$ is finite a.s.
  
Assume that $C(o)$ is infinite. Let $v$ be a vertex such that $C(o)
\subseteq S_v$. There are finitely many (no more than
$(d^{\ell(o)-\ell(v)+1}-1)/(d-1)$) vertices $u_i$ in $S_v$ such that
$\ell(u_i)\leq\ell(o)$, and the downwards cluster $C^+(u_i)$ of each such
vertex is finite a.s. On the other hand, $C(o)$ can be infinite only if
$o$ is connected to vertices $w$ on arbitrarily deep levels in $S_v$. If
we consider an open path from $o$ to such a $w$, then the first vertex on
this path which is on the level of $w$ is actually an element of
$C^+(u_i)$ for one of the vertices $u_i$. But this is impossible if $w$ is
located deep enough, hence $C(o)$ must be finite. \end{proof}

\chap{Diestel-Leader graphs}
\label{chap:43}

\sect{Definition}
\label{sec:43:defn}

\newcommand{\G}{\Gamma_{\alpha,\beta}}

Diestel and Leader \cite{DiLe:non-cayley} give the following example
of a graph with non-unimodular automorphism group.  They conjecture
that this transitive graph is not quasi-isometric to any Cayley graph.

Fix integers $\alpha, \beta \ge 2$.  Let $\tree{T}_\alpha$ and
$\tree{T}_\beta$ be an $(\alpha+1)$-regular and a $(\beta+1)$-regular tree,
respectively.  Choose an end of $\tree{T}_\alpha$ and an end of
$\tree{T}_\beta$, and orient the edges of each tree towards its distinguished
end.  Now construct the graph $\Gamma'_{\alpha,\beta} =
(V'_{\alpha,\beta},E'_{\alpha,\beta})$ with vertices
$V'_{\alpha,\beta}=V(\tree{T}_\alpha) \times V(\tree{T}_\beta)$ and edges
\[
\begin{split}
  E'_{\alpha,\beta} = \bigl\{ \{(u_1,u_2), (v_1,v_2)\} \bigm |
  &\text{$\{u_1,v_1\}\in E(\tree{T}_\alpha)$ and
    $\{u_2,v_2\}\in E(\tree{T}_\beta)$} \\
    &\text{are oriented in opposite directions.}
  \bigr\}.
\end{split}
\]

Note that if $\ell_1 : \tree{T}_\alpha\to \Ints$ and $\ell_2 :
\tree{T}_\beta\to \Ints$ are level functions for $\tree{T}_\alpha$ and
$\tree{T}_\beta$ respectively, then $\ell_1(u_1)+\ell_2(u_2)=
\ell_1(v_1)+\ell_2(v_2)$ for any edge $\{(u_1,u_2), (v_1,v_2)\}$ of 
$\Gamma'_{\alpha,\beta}$.  Thus, $\Gamma'_{\alpha,\beta}$ has 
infinitely many connected components, all isomorphic.  Define $\G$ to 
be one such connected component.

Figure~\ref{fig:43} illustrates a portion of~$\G$ when $\alpha=3$ and
$\beta=2$, along with a path in it.

\begin{figure}[ht]
  \begin{center}
    \includegraphics{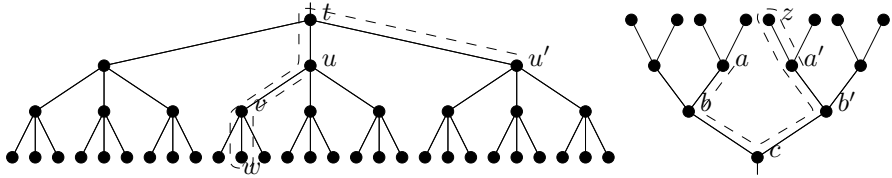}
    \caption{The Diestel-Leader graph $\Gamma_{3,2}$}
    \label{fig:43}
  \end{center}
  \begin{small}
    Nodes in $\G$ are pairs of vertices at \emph{the same} level; edges
    must follow both trees' edges.  Sample path: $(u,a), (v,b), (w,c),
    (v,b'), (u,a'), (t,z), (u',a')$.
  \end{small}
\end{figure}

We also define the two projections onto the first and
second components of $V_{\alpha,\beta}$, labelled $\pi_1: \G\to
\tree{T}_\alpha$ and $\pi_2: \G\to\tree{T}_\beta$, and note that if
$\{x,y\}$ is an edge of $\G$, then $\{\pi_1(x),\pi_1(y)\}$ and
$\{\pi_2(x),\pi_2(y)\}$ are edges of $T_\alpha$ and $T_\beta$
respectively.  Also, define a level function $\ell \bydef \ell_1 \circ
\pi_1$. The level sets $L_k$ for $k\in \Ints$ are defined relative to this 
function $\ell$.
We shall refer to an edge in $\G$ from $x$ to $y$ as going \defn{up} if
$\ell(y) = \ell(x) - 1$, and going \defn{down} if $\ell(y) = \ell(x) + 1$.

Here are some standard facts concerning~$\G$, which appear in
\cite{DiLe:non-cayley} and \cite{Woess:book}.

\begin{itemize}
\item $\G$ is clearly a transitive graph.
  
\item $\aut{\G}$ is unimodular iff $\alpha = \beta$, as the
  combinatorial characterization is easily checked. 
  
\item $\aut{\G}$ is always amenable. This is because $\aut{\G}$ is a
subgroup of the
  direct product of the groups of those automorphims of $T_\alpha$ and
  $T_\beta$ that preserve the distinguished end. As we have seen, these two
  groups are amenable, and group amenability is preserved by direct
  sums and by going to a subgroup.

\item $\G$ is amenable iff $\alpha= \beta$. This follows from the
  facts above and Lemma \ref{lem:amenunim}.
\end{itemize}

We will prove the following:

\begin{theorem}\label{thm:43:dead}
  If $\alpha > \beta$, then critical percolation on $\G$ almost surely has
no infinite
  clusters. Furthermore, $\Probp{p_c}{x\longleftrightarrow y}\leq
C\rho^{\dist (x,y)}$ for any $\rho\in \left(\frac{\beta}{\alpha},1\right)$
and suitable $C=C(\rho)<\infty$. \end{theorem}

 Many of the lemmas used in the proof are primarily
 combinatorial, and hold also when $\alpha=\beta$.  So we shall not assume 
 $\alpha\neq\beta$ unless it is explicitly stated in a lemma. The
unimodular case $\alpha=\beta$ will be discussed in
Section~\ref{chap:lamp}. 
The main help in the proofs will be that the
geometry of the graph $\G$ has some similarities with that of a tree:

\begin{note}\label{note:simple_path}
  Let $x_0,x_1,\ldots$ be a path in $\G$, such that the edge from
  $x_0$ to $x_1$ goes down, and  $\forall i\geq 0: \ell(x_i)\geq \ell(x_0)$. 
  Then the path $\pi_1(x_0), \pi_1(x_1), \ldots$ stays within the
  downwards subtree of $\pi_1(x_0)$.
\end{note}

This motivates the following definitions:

\begin{definition}\label{def:downforward}
 The \defn{forward subcluster} of a vertex $o\in V(\G)$ is
 the set  
\[\begin{split}
C^+(o)=C^+(o,\omega) \bydef \{ v \st o \longleftrightarrow
 v &\text{ by an open path } (o=v_0,v_1,\dots,v_n=v)\\
&\text{ such that }\ell(o)\leq\ell(v_i)\leq\ell(v)\text{ for all }i\}.
\end{split}\]
 Furthermore, the \defn{downwards subcluster} of $o$ is
  \[\begin{split}
C'(o) = C'(o,\omega) \bydef \{ v \st o \longleftrightarrow
 v &\text{ by an open path } (o=v_0,v_1,\dots,v_n=v)\\
&\text{ such that }\ell(o)\leq\ell(v_i)\text{ for all }i\}.
\end{split}\]
\end{definition}

\sect{Finiteness downwards}
\label{sec:43:down}

As before, the first step of the proof is to bound the rate of growth of
the forward part of the critical cluster, and to conclude that the cluster
cannot directly go down infinitely deeply.

\begin{lemma}\label{lem:43:bounded}
  Let $o\in V(\G)$ and consider critical percolation on $\G$.  Define
  $V_k^+ = V_k^+(o) \bydef \left\{ x : x\in C^+(o), \ell(x)-\ell(o) = k
  \right\}$.  Then the values $e_k = e_k(p_c) \bydef
  \Expectation{\sizeof{V_k^+}}$ satisfy $e_k \le 1$.
\end{lemma}

\begin{proof} We can copy the proof of Lemma \ref{lem:dt:bb}.  The only
difference is that the independence of the number of offspring of any two
vertices on the same level of the Galton-Watson process we are building is
now provided by Note \ref{note:simple_path}, as opposed to the earlier
explicit restriction that the paths in $C^+(o)$ should stay inside the
subtree $S_o$. 
\end{proof}

\begin{lemma}\label{lem:43:finite_proj}
  In critical percolation on $\G$, we have that $\pi_1[C'(o)]$ is finite a.s.
\end{lemma}

\begin{proof} Exactly as in the proof of Lemma \ref{lem:dt:posfinite}, we
can use Fatou's lemma and the sequence of events $\GG_k$ to conclude that
there is a random integer $K$ such that $V_k^+(o)$ is empty for all $k>K$.
In other words, $\pi_1(C^+(o))$ is finite almost surely. Now note that,
unlike $C^+(o)$, the set $C'(o)$ is necessarily connected.  Hence, if
$\pi_1(C'(o))$ was infinite, then for any $k>0$, there would be a simple
open path in $C'(o)$ between $o$ and some vertex of $L_k$.  The first time
this path enters $L_k$, at vertex, say, $v\in L_k$, then $v\in C^+(o)$
would also hold.  Since $k$ was arbitrary, $\pi_1(C^+(o))$ would be
infinite, too. 
\end{proof}

The last result can be strengthened by the following simple lemma, which
shows that the structure of $\G$ ensures that a connected set $C$ with
finite $\pi_1(C)$ cannot be infinite.

\begin{lemma}\label{lem:43:finite_vertices}
  Let~$C$ be a connected component of~$\G$ such that there exists
  some~$k\in\Ints$ with $\ell(q) \le k$ for any $q\in C$.  Then for all
  $u\in\tree{T}_\alpha$, the projection $\pi_1$ maps only finitely many
elements of~$C$ to~$u$, and indeed
  $\sizeof{\pi_1^{-1}[u]\cap C} \le \beta^{k-\ell_1(u)}$.
\end{lemma}

\begin{proof} If $\ell_1(u) > k$, then $\pi_1^{-1}[u]\cap C = \emptyset$.  
Now suppose that this set is non-empty, with $k-\ell_1(u)=j \geq 0$, and
take some $(u,a) \in \pi_1^{-1}[u]\cap C$, with $a\in T_\beta$. Consider
the unique ancestor $b\in T_\beta$ of $a$ which has
$\ell_2(b)=\ell_2(a)-j$, and let $T\subset T_\beta$ be the infinite
subtree of descendants of $b$. Denote the $j$th descendants of $b$ by
$a_1=a,a_2,\dots,a_{\beta^j} \in T$.  Note that any open path
$\gamma\subseteq C$ starting from $(u,a)$, because of $\ell(q)\le k$ for
all $q\in\gamma$, satisfies $\pi_2(\gamma)\subseteq T$.  Hence
$\pi_1^{-1}[u]\cap C\subseteq \{(u,a_1),\dots,(u,a_{\beta^j})\}$, and the
claim follows. \end{proof}

\sect{Finiteness upwards}
\label{sec:43:up}

Next we prove that almost surely no vertex can connect to vertices
unboundedly ``upwards'' of it in the tree.

\begin{lemma}\label{lem:43:upward_bound}
  Consider critical percolation on $\G$, and fix a vertex~$o$.  For all
$k\in\Nats$, define 
  \[
  V^-_k(o) \bydef \set{x\st \ell(o)-\ell(x)=k, \ o\in C^+(x)}.
  \]
  Then $\Expectation\sizeof{V^-_k} \le \left(\frac{\beta}{\alpha}\right)^k$.
\end{lemma}

\begin{proof}
 There are $\beta^k$ vertices in $\G$ that have a positive probability to
appear in $V^-_k$, and all these probabilities are the same, also equaling
to $q_i\bydef\Probp{p_c}{\exists u \in C^+(o)
  \text{ with }\pi_1(u)=a_i}$, where the $a_i$, $i=1,\dots,\alpha^k$, are
the $k$'th generation descendants of $\pi_1(o)$ in $T_\alpha$. Now,
rewriting $e_k$ from Lemma \ref{lem:43:bounded} as
  \[
  1 \ge e_k = \sum_{i=1}^{\alpha^k} q_i = \alpha^k q_1
  \] gives $q_i\leq \alpha^{-k}$, and the desired bound follows from the
linearity of expectation.  
\end{proof}

\begin{lemma}\label{lem:43:finiteness}
  Suppose $\alpha > \beta$, and for all $k\in\Nats$, define
 \[
  U^-_k(o) = \set{x\st \ell(o)-\ell(x)=k, \ o\in C'(x)}.
  \]
Then $a_k\bydef \Expectation\sizeof{U^-_k}<\infty$. 
 \end{lemma}

\begin{proof} First we prove that $a_0\leq
\frac{\alpha}{\alpha-\beta}<\infty$, then that $a_{k+1}\leq a_k
(\beta/\alpha) a_0<\infty$ for all $k\in\Nats$.

Consider a simple open path $\gamma$ connecting $o$ to $x$ and showing
$x\in U^-_0$.  Let~$y$ be the \emph{last} lowest vertex on the path.  
Write $j = \ell(y)-\ell(o)$.  Then the portion of $\gamma$ between $y$ and
$x$ shows that $x\in V^-_{j}(y)$, with the definition of Lemma
\ref{lem:43:upward_bound}, while the portion of $\gamma$ between $y$ and
$o$ shows that $y\in C^+(o)$. Now, such an open path $\gamma$, going
through these vertices $o,y,x$, though with $y$ not being necessarily the
\emph{last} lowest vertex, exists if and only if both events $\{x\in
V^-_{j}(y)\}$ and $\{y\in C^+(o)\}$ occur, due to disjoint sets of open
edges; i.e.~if{f} $\{x\in V^-_{j}(y)\} \,\square\, \{y\in C^+(o)\}$
happens, with the notation of the van den Berg -- Kesten inequality, see
\cite{BK:1985} or \cite{grimmbook}.

This BK inequality says that for increasing measurable events $\A$ and
$\B$ in independent $p$-percolation,
$\Probp{p}{\A\square\B}\leq\Probp{p}{\A}\Probp{p}{\B}$. Therefore, at
$p_c$, 
\[ 
\begin{split} 
a_0=\Expectation{|U^-_0|} &\leq
\sum_{j=0}^{\infty} \sum_x \sum_y \Prob{\{x\in V^-_{j}(y)\} \,\square\,
\{y\in C^+(o)\}} \\ &\leq\sum_{j=0}^\infty\sum_x \sum_y \Prob{x\in
V^-_{j}(y)}\Prob{y\in C^+(o)}\\ &\leq\sum_{j=0}^\infty
\left(\frac{\beta}{\alpha}\right)^{j}\cdot
1=\frac{\alpha}{\alpha-\beta}<\infty, 
\end{split} 
\] 
where the sums are
over $\{x:\ell(x)=\ell(o)\}$ and $\{y:\ell(y)=\ell(o)+j\}$, and we used
Lemmas \ref{lem:43:upward_bound} and \ref{lem:43:bounded} to get the third
line.

Now take a simple open path $\gamma$ from $o$ to $x$ and showing $x\in
U^-_{k+1}(o)$.  Let~$z$ be the \emph{first} vertex on this path that lies
in $U^-_{k+1}(o)$, and let $t$ be previous vertex on the path. Then the
portion of $\gamma$ between $o$ and $t$ shows that $t\in U^-_{k}(o)$,
while the portion of $\gamma$ between $z$ and $x$ shows that $x\in
U^-_0(z)$.
 Similarly as above, such an open path $\gamma$, going through these
vertices $o,t,z,x$, exists if and only if the three events $\{t\in
U^-_{k}(o)\}$, $\{(t,z)\text{ is open}\}$ and $\{x\in U^-_{0}(z)\}$ occur
on disjoint sets of open edges. Hence the BK inequality now gives 
\[
a_{k+1}\leq a_k \cdot (\beta p_c) \cdot a_0. 
\] 
Since $T_\alpha$ is a subgraph of $\G$, we have $p_c\leq 1/\alpha$. 
Therefore, by induction, $a_k\leq (\beta/\alpha)^k a_0^{k+1}\leq
\left(\frac{\beta}{\alpha-\beta}\right)^k\frac{\alpha}{\alpha-\beta}<\infty$.
\end{proof}

Now plugging the finiteness of the $a_k$'s into a similar, but more
refined argument, we get for all $\alpha>\beta$ that $a_k\to 0$
exponentially, as $k\to\infty$.

\begin{lemma}\label{lem:43:top_bound}
  Suppose $\alpha > \beta$.  Define the event
  \[
  \A_k = \set{\text{$o$ is
      connected by an open path to a vertex~$k$ levels above it}}.
  \]
   Then, in critical percolation, $\lim_{k\to\infty} \Prob{\A_k} = 0$, 
decaying exponentially.
\end{lemma}

\begin{proof}
  Note that $\A_k = \set{U^-_k \neq \emptyset}$.  Thus, by Markov's
inequality, it is enough to show that for $a_k\bydef
\Expectation\sizeof{U^-_k}$ we have $\lim_{k\to\infty} a_k=0$
exponentially quickly.

 Consider a simple open path $\gamma$ connecting $o$ to $x$ and showing
$x\in U^-_k$.
  Let~$y$ be the \emph{last} lowest vertex on the path.  Write $j
  = \ell(y)-\ell(o)$.  Then the portion of $\gamma$ between $y$ and $x$
shows that $x\in V^-_{j+k}(y)$, with the definition of Lemma
\ref{lem:43:upward_bound}.
 Now let $z$ be the \emph{last} highest vertex on the portion of $\gamma$
between
  $o$ and $y$, and write $i = \ell(o)-\ell(z)$. Clearly, $0\le i\leq k$,
and $z\in U^-_i(o)$.  The path $\gamma$ also shows that $y\in C^+(z)$.

The existence of such an open path $\gamma$, going through these vertices
$x,y,z$, is equivalent to the occurrence of the three events
 $\{x\in V^-_{j+k}(y)\}$, $\{y\in C^+(z)\}$ and $\{z\in U^-_i(o)\}$ on
disjoint edge sets. Hence the BK inequality gives that, at $p_c$,
\[
\begin{split} 
\Expectation{|U^-_k|} &\leq \sum_{i=0}^k \sum_{j=0}^{\infty}
\sum_x \sum_y \sum_z \Prob{\{x\in V^-_{j+k}(y)\} \,\square\, \{y\in
C^+(z)\} \,\square\,\{z\in U^-_i(o)\}}\\
&\leq\sum_{i=0}^k\sum_{j=0}^\infty\sum_x \sum_y \sum_z \Prob{x\in
V^-_{j+k}(y)}\Prob{y\in C^+(z)}\Prob{z\in U^-_i(o)}\\
&\leq\sum_{i=0}^k\sum_{j=0}^\infty a_i\cdot 1\cdot
\left(\frac{\beta}{\alpha}\right)^{j+k}
=\frac{\alpha}{\alpha-\beta}\left(\frac{\beta}{\alpha}\right)^k\sum_{i=0}^k
a_i, 
\end{split} 
\] 
where the sums are over $\{x:\ell(x)=\ell(o)-k\}$,
$\{y:\ell(y)=\ell(o)+j\}$, and $\{z:\ell(z)=\ell(o)-i\}$, and we used the
definition of $a_i$ and Lemmas \ref{lem:43:upward_bound} and
 \ref{lem:43:bounded} to get the third line.

The finiteness of the $a_i$'s is known from Lemma \ref{lem:43:finiteness}.
Now suppose $\limsup_{i\to\infty} a_i\geq \delta>0$. Because of the
exponential decay of the factor $\left(\frac{\beta}{\alpha}\right)^k$ in
the previous inequality, this can happen only if $\limsup_{i\to\infty}
a_i=\infty$. But then there are infinitely many indices $m$ for which
$a_m=\max\{a_0,a_1,\dots, a_m\}$, and for such an $m$ our inequality
implies $a_m\leq\frac{\alpha}{\alpha-\beta} (m+1)
a_m\left(\frac{\beta}{\alpha}\right)^m$. But this is impossible if $m$ is
large enough.  Hence $\lim_{k\to\infty} a_k=0$. Moreover, convergence
implies boundedness, $a_k\leq A$ for some $A$, hence we actually have
$a_k\leq
\frac{\alpha}{\alpha-\beta}(k+1)A\left(\frac{\beta}{\alpha}\right)^k$,
which is less than $B\rho^k$ for $\frac{\beta}{\alpha}<\rho<1$ and $B>0$
large enough. \end{proof}

Iterating further our argument gives the following:
\begin{lemma}\label{lem:43:total} Suppose $\alpha > \beta$, and for all
$k\in\Nats$, define \[ W^-_k(o):=L_{-k}\cap C(o). \] Then, in critical
percolation, $\Expectation{\sizeof{W^-_k(o)}}<\infty$, and they decay
exponentially in $k$. \end{lemma}

\begin{proof} Consider an open path $\gamma$ from $o$ to $x\in W^-_k$. Let
$w_1$ be the last highest vertex on $\gamma$, and $v_1$ be the last lowest
vertex on the portion of $\gamma$ from $w_1$ to $x$. Then, for $t\geq 2$,
let $w_t$ be the last highest vertex on the portion of $\gamma$ from
$v_{t-1}$ to $x$, and $v_t$ be the last lowest vertex on the portion of
$\gamma$ from $w_t$ to $x$. We make these definitions for all $t\geq 1$,
but there will certainly be a smallest $T\geq 1$ such that $w_t=v_t=x$ for
all $t\geq T$. Writing $j_t=\ell(x)-\ell(w_t)$ and
$i_t=\ell(v_t)-\ell(x)$, we have $j_1 > j_2 > \dots > j_{T-1} >
j_T=j_{T+1}=\dots=0$ and $i_1 > i_2 > \dots > i_{T-1} \geq
i_T=i_{T+1}=\dots =0$.  Note that $w_1\in U_{k+j_1}^-(o)$ and $w_{t+1} \in
V_{i_t+j_{t+1}}^-(v_t)$, while $v_t \in C^+(w_t)$.  The BK inequality now
gives 
\[ 
\begin{split} 
\Expectation{|W^-_k|} & \leq \sum_{j_1>j_2>\dots\atop i_1>i_2>\dots}
\sum_{w_1,w_2,\dots \atop v_1,v_2,\dots}\, 
\Prob{w_1\in U^-_{k+j_1}(o)}\prod_{t=1}^\infty \Bigl(\Prob{w_{t+1} \in
V_{i_t+j_{t+1}}^-(v_t)}\Prob{v_t\in C^+(w_t)}\Bigr)\\ 
&\leq \sum_{j_1>j_2>\dots \atop i_1>i_2>\dots}
 B \rho^{k+j_1} \prod_{t=1}^\infty
\left(\frac{\beta}{\alpha}\right)^{i_t+j_{t+1}}\\ 
&=B \rho^k \frac{1}{1-\rho}\Biggl( \sum_{j_2>j_3>\dots}
\left(\frac{\beta}{\alpha}\right)^{j_2+j_3+\dots}\Biggr) 
\Biggl(\sum_{i_1>i_2>\dots}
\left(\frac{\beta}{\alpha}\right)^{i_1+i_2+\dots}\Biggr)\\
&= B \rho^k \frac{1}{1-\rho} \left( \sum_{n=0}^\infty
\left(\frac{\beta}{\alpha}\right)^n q(n)\right)^2, 
\end{split} 
\] 
where, to get the second line, we again used Lemmas
\ref{lem:43:upward_bound} and 
\ref{lem:43:bounded} and wrote $a_k\leq B\rho^k$ from the proof of
Lemma~\ref{lem:43:top_bound}, while, in the last line, we wrote $q(n)$ for
the number of partitions of $n\in\Nats$ with all distinct parts, with the
convention $q(0)=1$. It is easy to see that $q(n)$ has subexponential
growth, 
$$q(n)\leq \sum_{k=1}^{\sqrt{2n}} {n\choose k} \leq 
\exp(C\sqrt{n}\log n),$$ 
but very precise estimates exist: it is well-known \cite{Andrews} that $q(n)$
is also the number of partitions of $n$ into odd parts, and we have
$$q(n)\sim \frac{e^{\pi\sqrt{n/3}}}{4\cdot 3^{1/4}n^{3/4}},$$ see
\cite{Iseki, Hagis}. We thus conclude that the last infinite series
converges to a finite value $Q_{\alpha,\beta}$ for any $\beta<\alpha$.
That is, $$\Expectation{|W^-_k|}\leq
B\rho^k\frac{1}{1-\rho}Q_{\alpha,\beta}^2,$$ and the proof is complete.
\end{proof}

\begin{proof}[Proof of Theorem~\ref{thm:43:dead}]
  Consider a component $C = C(o)$ of critical percolation on $\G$.  In
  view of Lemma \ref{lem:43:finite_vertices}, it suffices to show
  that a.s.~$\pi_1(C)$ is finite to conclude that a.s.~$C$ is finite.
  
 By Lemma~\ref{lem:43:top_bound}, a.s.~every component~$C$ has a
  highest level, which contains a finite number of vertices.  By
  Lemmas \ref{lem:43:finite_proj} and \ref{lem:43:finite_vertices},
  a.s.~each vertex has a finite downwards subcluster.  But $C$ is just 
  the union of the downwards subclusters of its vertices at the highest
  level, hence is (a.s.) finite.

The exponential decay of the connection probabilities follows
immediately from Lemma \ref{lem:43:total} and the fact that as a function of
$t\in\Nats$, there is an exponentially large number of vertices $v$ with
the properties that $\dist(o,v)=t$, all $v$'s are on the same level of
$\Gamma_{\alpha,\beta}$, and, moreover, their connection probabilities
to $o$ are the same.
\end{proof}

The characterization of $p_u$ due to Schonmann~\cite{Schon:stab}
and the amenability of $\Gamma_{\beta,\beta} $ imply that 
$$
p_u(\G) \le p_u(\Gamma_{\beta,\beta}) =  
p_c(\Gamma_{\beta,\beta})\le 1/\beta
$$
for $\beta \leq \alpha$.
A condition on $\alpha$ and
$\beta$ for $p_c(\G)<p_u(\G)$ can be easily given using a result of
Schramm, see \cite[Theorem 6.28]{LPbook}. If we denote by $a_n(G)$
the number of simple loops of length $n$ containing a fixed vertex
$o\in V(G)$, and $\gamma(G):=\limsup_n a_n(G)^{1/n}$, then $p_u(G)\geq
1/\gamma(G)$.  For $G=\G$ it is not difficult to see that 
$$\gamma \leq  \sqrt{\alpha\beta}+\sqrt{(\alpha-1)(\beta-1)},$$
by the following argument.

First of all, $\G$ is a bipartite graph, so $a_{2n+1}(\G)=0$, while
$a_{2n}(\G)$ is bounded from above by the number of simple
non-backtracking paths of length $2n$ ending on the starting level.  
(Note that in this estimate we do not lose much by relaxing the
loop-condition; however, excluding immediate backtracks is quite far from
ensuring that the path be simple.) In such a path, we have $n$ upwards and
$n$ downwards moves, in an arbitrary order, with $k$ instances of changing
direction from upwards to downwards, where $k\in\{0,1,\dots,n\}$. Then,
the number of changes in direction from downwards to upwards is between
$k-1$ and $k+1$. The number of such sequences with a given $k$ value is at
most $Cn{2n\choose 2k}$. When such a path changes direction, to avoid
backtracking, it has $\alpha-1$ or $\beta-1$ ways to continue; when it
does not change direction, it has $\alpha$ or $\beta$ ways. Therefore, 
\[
\begin{split} a_{2n}(\G) & \leq \sum_{k=0}^n 
C'n{2n\choose 2k}\alpha^{n-k}\beta^{n-k}(\alpha-1)^k(\beta-1)^k\\ 
& < C'n(\alpha\beta)^n \sum_{k=0}^{2n} {2n\choose k} x^k, 
\qquad\hbox{with}\quad x=\sqrt{\frac{(\alpha-1)(\beta-1)}{\alpha\beta}},\\ 
& = C'n\left(\sqrt{\alpha\beta}\right)^{2n}(1+x)^{2n}. 
\end{split} 
\] 
Taking the $(2n)$th root of the last line gives the claimed bound on
$\gamma$. 

On the other hand, it is clear that $p_c\leq 1/\alpha$.
(By considering small cycles, this inequality, as well as the above bound on
$\gamma$, can be improved.) Hence 
$$\sqrt{\alpha\beta}+\sqrt{(\alpha-1)(\beta-1)}\leq \alpha 
\qquad\hbox{implies}\qquad p_c(\G)<p_u(\G).$$ 
This is the case e.g.~for $\Gamma_{6,2}$, and for $\alpha\geq 4\beta$, 
in general.

If one could deduce from the uniform exponential decay of connection
probabilities at $p_c$ (which we have verified for all $\alpha<\beta$)
that for some $p>p_c$, the connection
probabilities still tend to 0, it would follow that  $p_c(\G)<p_u(\G)$  by the
Harris-FKG inequality. 

\chap{The lamplighter group}
\label{chap:lamp}

\newcommand{\gamhalf}{\Gamma^+_{\alpha,\alpha}}

Recall that when $\alpha = \beta$, the graph $\G$ is
 amenable and unimodular.  The first half of our proof of Theorem 
\ref{thm:43:dead} still holds,  
 but the bound of Lemma \ref{lem:43:upward_bound} does not mean 
exponential decay, 
and so this method brakes down. 

Take the ``positive part'' $\gamhalf$ defined by taking the subgraph
induced by the vertices
\[
V(\gamhalf) = \set{v\in V(\Gamma_{\alpha,\alpha})\st
  \text{$\pi_1(v)$ is a descendant of $\pi_1(o)$}}.
\]
 Clearly, $p_c(\Gamma_{\alpha,\alpha})\leq p_c(\gamhalf)$, and our proof
above shows that $p_c(\gamhalf)$-percolation on $\gamhalf$ has no infinite
clusters. This remains true for $p_c(\Gamma_{\alpha,\alpha})$-percolation
on $\gamhalf$, so any infinite path in
$p_c(\Gamma_{\alpha,\alpha})$-percolation on $\Gamma_{\alpha,\alpha}$
would have to cross the plane $\set{v\st \ell(v)=\ell(o)}$ infinitely many
times.

A special interest in the graphs $\Gamma_{\alpha,\alpha}$ comes from the
fact that they also arise as Cayley graphs of the so-called ``lamplighter
groups'', introduced by Ka{\u\i}manovich and Vershik (Example 6.1 of
\cite{KV:lamplighter:1983}), and further studied from a probabilistic
point of view e.g.~by Lyons, Pemantle and Peres
\cite{LPP:lamplighter:1996} and Woess \cite{Woess:DL}.

\begin{definition}[Example~6.1 of \cite{KV:lamplighter:1983}] Consider the
direct sum $\sum_{\Ints^k} \Ints_2$, which can also be viewed as the
additive group $F_0(\Ints^k,\Ints_2)$ of finitely supported
$\{0,1\}$-configurations on $\Ints^k$, with
  the operation of pointwise addition {\rm mod} 2.  The value of a
   configuration $f\in F_0(\Ints^k,\Ints_2)$ on an element $x\in\Ints^k$ 
 will be denoted by $f(x)$ and the \defn{support} $\set{x\in \Ints^k: f(x) 
 \neq 0}$ of $f$ by $\supp f$.
  Let
  \[
  G_k = \Ints^k \ltimes F_0(\Ints^k,\Ints_2)
  \]
  be the
  \emph{semidirect} product of the groups $\Ints^k$ and
  $F_0(\Ints^k,\Ints_2)$, where the lattice $\Ints^k$ acts naturally on
  $F_0(\Ints^k, \Ints_2)$ by shifts.
\end{definition}

\begin{figure}[htbp]
  \begin{center}
    \includegraphics{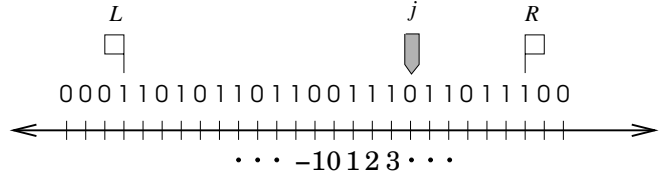}
    \caption{A configuration of lamps and the lamplighter in $G_1$.}
    \label{fig:lampconfig}
  \end{center}
\end{figure}

The group $G_1$ was named the \defn{lamplighter group} because of the 
following interpretation. 
Imagine a lamplighter standing on an infinite street, with
lamps at every integer coordinate.  Any element $(j, f)$ describes a
configuration: the lamplighter is next to lamp~$j$, and~$f$ is the
indicator function of the finite set~$F$ of lamps which are lit.  For
convenience, we shall also denote this element by $(j,F)$. Define the 
 \defn{left} and \defn{right flag} functions by $L((j,F))\bydef\min F$ and 
 $R((j,F))\bydef \max F$,   with $\min \{\}\bydef +\infty$, $\max 
 \{\}\bydef -\infty$ for the empty set $\{\}$, and the \defn{lamplighter 
position} by $\ell((j,F))\bydef j$. See Figure~\ref{fig:lampconfig}.

The group operation is given by $(j,F)\cdot (j',F') = (j+j', F \symdiff
(j+F'))$, where $\symdiff$ is symmetric set difference: the
lamplighter flips the lamps~$F'$ relative to her current position, and
advances $j'$ lamps.  

Recall the construction of the graph~$\Gamma_{2,2}$ by orienting two
$3$-regular trees in opposite directions.  Label the edges of each
tree `0' or `1', to satisfy these conditions:
\begin{enumerate}
\item The two ``downwards'' edges from each vertex are labelled `0'
  and `1';
\item The edges of every ``upwards'' path $v_0, v_1, \ldots$ are
  eventually all labelled `0'.
\end{enumerate}
Then, given any vertex~$v$ at level $\ell(v)=k$, we may identify~$v$
with the element $(k,f)$ of~$G_1$ as follows: Let $a_0,a_1,\ldots$ and
$b_0,b_1,\ldots$ be the labels of the edges along the paths
upwards from $\pi_1(v)$ and $\pi_2(v)$, respectively.  For $j\ge 0$,
define $f(k+j)=b_j$ and $f(k-1-j)=a_j$.  Then~$f$ has finite support,
so $(k,f)$ is in $G_1$ indeed.  In fact, $\Gamma_{2,2}$ is the Cayley graph of
the lamplighter group $G_1$ with generators $\set{(\pm 1,\set{0}), (\pm 
1,\{\})}$.

\begin{figure}[htbp]
  \begin{center}
    \includegraphics{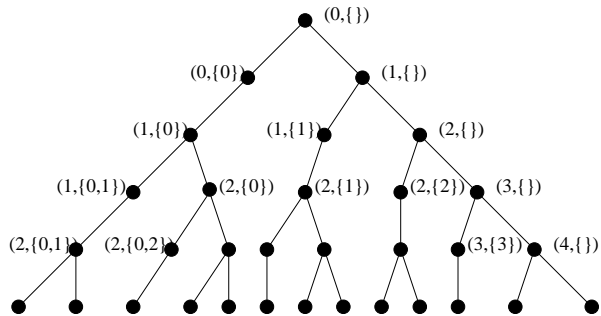}
    \caption{The Fibonacci tree in the lamplighter group.}
    \label{fig:fibonacci}
  \end{center}
\end{figure}

Another natural Cayley graph~$G$ is given by the generators 
$(0,\set{0})$ (the lamplighter flips the state
of the current lamp and stays in place) and $(\pm 1,\{\})$ (the
lamplighter advances one lamp). Consider again the ``positive half'' 
$G^{+}$ of~$G$, 
defined by taking only the vertices
\[
\set{(j,F)\st j\ge 0, \forall k\in F: k\ge 0}.
\]
This is the portion of the graph accessible to the lamplighter if she
 is limited to the non-negative portion of the street. The subset of 
$V(G^+)$ given by $\{v: R(v)\leq \ell(v)\}$ induces a tree $\FF$, the 
 so-called Fibonacci tree, identified in \cite{LPP:lamplighter:1996}. See 
Figure~\ref{fig:fibonacci}.

We again have $p_c(G)\leq p_c(G^+)$, and will consider
$p_c(G^+)$-percolation on $G^+$. For a vertex $o\in V(\FF)$, define the
\defn{forward cluster} $C^+(o)=C^+(o,\omega)$ as the set of vertices $v\in
V(\FF)$ accessible by open paths $o=v_0,v_1,\dots,v_n=v$ inside $G^+$ (not
necessarily inside $\FF$) in a $p_c(G^+)$-percolation configuration
$\omega$, such that $\ell(o)\leq \ell(v_i) \leq \ell(v)$, and the lamp at
$R(o)$ is never adjusted in the path.

It is easy to see that we have the required independence in order to make
our usual Galton-Watson argument work, therefore $e_k(p_c)\bydef
\Expectationp{p_c}{|C^+(o) \cap \{v:\ell(v)=k\}|} \leq 1$ for all $k\geq
0$. Again, as in Lemmas \ref{lem:dt:bb} and \ref{lem:43:finite_proj}, we
can conclude that $C^+(o)$ must be finite. Moreover, any open infinite
simple path from $o$ in $G^+$ would have infinitely many vertices inside
$C^+(o)$, therefore the whole cluster of $o$ is almost surely finite in
critical percolation on $G^+$.

We have shown two transitive amenable graphs for which we know that
critical percolation on the ``positive part'' almost surely has no
infinite clusters, but we cannot prove this for the whole graph.
Analogously, Barsky, Grimmett and Newman \cite{BGN:half-space:1991} proved
that critical percolation on the half-space graphs of the integer lattices
$\Ints^d$ has no infinite clusters.

  \chap{Acknowledgments}
  \label{chap:thanks}
  
The starting point of this paper was Ariel Scolnicov's Master's thesis,
written in 2000 at the Department of Mathematics at the Hebrew University,
under the supervision of Yuval Peres.  We thank Itai Benjamini, Noam
Berger, Elchanan Mossel, Asaf Nachmias, \'Ad\'am Tim\'ar and Tamar Ziegler
for useful comments, and the referee for many important corrections to an
earlier version of the paper.

\bibliographystyle{alpha}
\bibliography{papers}
\vskip 0.5 cm

\end{document}